\documentclass{article}

\def\PSL{{\rm PSL}}
\def\F{\hbox{\bf F}}
\def\Q{\hbox{\bf Q}}
\def\Z{\hbox{\bf Z}}
\def\C{\hbox{\bf C}}
\def\Card{{\rm Card}\;}
\def\PGL{\mathop{\rm PGL}\nolimits}

\def\L{\mathop{\rm L}\nolimits}

\def\ni{\noindent}
\def\ms{\mapsto} 
\def\ra{\rightarrow}
\def\eps{\varepsilon}
\def\fracb#1#2{\frac{{\textstyle #1}}{{\textstyle #2}}}
\def\mod{\mathop{\;\rm mod}\nolimits}

\begin{document}
\begin{center}
{\large {\bf Correspondances compatibles avec une relation binaire, rel\`evement d'extensions de groupe de Galois $\L_3(2)$ et probl\`eme de Noether pour $\L_3(2)$}}

\vspace{5ex}

{\sc Jean-Fran\c cois Mestre}
\vspace{5ex}

Universit\'e Paris 7,\\ 2 place Jussieu,
75230 Paris\\mestre@math.jussieu.fr 
\end{center}

\medskip

\section{Introduction}

Dans cet article, nous d\'efinissons la notion de correspondance compatible avec une relation binaire.
Dans le cas des relations d'incidence points-droites du plan projectif sur le corps fini $\F_2$, nous montrons
qu'il existe une telle correspondance "g\'en\'erique", et nous appliquons ce r\'esultat au probl\`eme de
Beckmann-Black et au probl\`eme de Noether relatif au groupe $\PGL_3(\F_2)$, not\'e parfois
$\L_3(2)$.

\medskip

{\sc D\'efinition.-} {\it Soit deux ensembles finis $A$ et $B$ de cardinaux respectifs $N$ et $M$,  $C$ une partie de $A\times B$,
$p_A:\;C\ra A$ et $p_B:\;C\ra B$ les projections. On suppose qu'il existe  deux entiers $n$ et $m$ tels
que, pour tout $a\in A$ et $b\in B$, $\Card p_A^{-1}(a)=n$ et $\Card p_B^{-1}(b)=m.$

Soit de plus $C_1$ et $C_2$ deux courbes d\'efinies sur un corps $k$. 
Une correspondance $\Gamma\subset C_1\times C_2$ de bidegr\'e $n-m$ est dite {\it compatible avec  $C$} si,
\'etant donn\'e un point g\'en\'erique $P$ de $C_1$, il existe une famille de $N$ points distincts ${\cal P}$ de $C_1$ 
contenant
$P$ et indic\'ee par $A$, ainsi qu'une famille de $M$ points distincts ${\cal Q}$ de $C_2$  indic\'ee par $B$, telle 
que $(a,b)\in C$ \'equivaut \`a  $(P_a,Q_b)\in \Gamma$.
}

\medskip

{\sc Exemple.-} Soit $n\geq 5$ un entier, et  $A$ (resp. $B$) les $n$ sommets (resp, les $n$ c\^ot\'es) d'un polygone convexe
\`a $n$ sommets, la partie $C\subset A\times B$ \'etant form\'ee des couples $(a,b)$ tels que
$a\in b$. Soient $C_1$ et $C_2$ deux coniques d\'efinies sur $\C$, et $\Gamma$ la correspondance
de bidegr\'e $2-2$
sur $C_1\times C_2$ d\'efinie par $(a,b)\in \Gamma$ si et seulement si $a$ appartient
\`a la tangente en $b$ \`a $C_2$. Supposons 
qu'il existe  une partie $M$ form\'ee de $n$ points distincts   de $\Gamma$  tels qu'on peut
indicer la projection de $M$ sur $C_1$ par $A$(resp. sur $C_2$ par $B$) de telle fa\c con que
$(a,b)\in C$ si et seulement si $(P_a,Q_b)\in M$. Le th\'eor\`eme de Poncelet affirme qu'alors
$\Gamma$ est compatible avec $C$.
 
\medskip

Nous montrons ici que, pour tout choix d'ind\'etermin\'ees $(x_a)$ indic\'ees par le plan
projectif $P^2(\F_2)$, il existe une correspondance sur $P^1\times P^1$ d\'efinie sur $\Q((x_a)_{a\in A})$, compatible
avec les relations d'incidence points-droites de $P^2(\F_2)$.

Plus pr\'ecis\'ement :

{\sc Th\'eor\`eme 1.-} {\it  Soit $A=P^2(\F_2)$ et $B$ l'ensemble des sept droites de $A$.
Pour tout syst\`eme d'ind\'etermin\'ees $(x_a)_{a\in A}$, il existe  une famille
$(y_b)_{b\in B}$ d'\'el\'ements distincts de $K=\Q((x_a))$, unique \`a homographies pr\`es, et
une  $3-3$ correspondance $\Gamma\subset P^1\times P^1$ d\'efinie sur $K$, unique une fois fix\'es
les $y_b$, 
telle que, pour tout $a\in A$ et $b\in B$, $(x_a,y_b)\in \Gamma$ si et seulement si $a\in b$, et  compatible avec les relations d'incidence points-droites
de $P^2(\F_2)$.
}

On en d\'eduit sans difficult\'e :

\medskip

{\sc Th\'eor\`eme 1'.-} {\it Il existe un polyn\^ome $S\in \Q[A_0,\ldots,A_6]$ tel que, si  $k$ est 
un corps de caract\'eristique nulle, si $P=X^7+a_6X^6+\ldots+a_0\in k[X]$ est un polyn\^ome
s\'eparable de groupe de Galois  contenu dans $\L_3(2)$ tel que $S(a_0,\ldots,a_7)\neq 0$,\footnote{dans ce genre de situation, nous dirons parfois ici pour abr\'eger
que {\it $P$ est  suffisamment g\'en\'eral}.} et 
si $f$ est une bijection de $A=P^2(\F_2)$ sur les racines de $P$ pr\'eservant sa structure
galoisienne (i.e. $a,b,c$ align\'es si et seulement si $f(c)\in k(f(a),f(b))$, il existe une famille
$(\beta_b)_{b\in B}$ d'\'el\'ements distincts du corps des racines de $P$, unique \`a homographies pr\`es,
et une correspondance $F\in k[X,Y]$, unique une fois fix\'es les $(\beta_b)$, v\'erifiant
$F(f(a),\beta_b)=0$ si et seulement si $a\in b$, et compatible avec les relations d'incidence
points-droites de $P^2(\F_2)$.
}

Nous donnons deux applications de ce th\'eor\`eme :

\bigskip

\begin{center}
I. {\sc Probl\`eme de Beckmann-Black}
\end{center}

\medskip

 Soit $k$ un corps, et $G$ un groupe fini. Le "probl\`eme inverse de Galois" pour le corps
$k$ et le groupe $G$ consiste
\`a rechercher s'il existe une extension galoisienne de $k$ 
de groupe de Galois $G$. Pour aborder ce probl\`eme, une m\'ethode
qui s'est r\'ev\'el\'ee f\'econde consiste \`a   
construire
une extension r\'eguli\`ere de $k(T)$ de groupe de Galois $G$.
Par sp\'ecialisation de $T$ en des \'el\'ement de $k$, on obtient
ainsi une infinit\'e d'extensions de $k$ de groupe de Galois $G$.

Il est naturel de se poser le probl\`eme inverse : soit $K/k$ une
extension de groupe de Galois $G$. Existe-il une extension galoisienne
$L/k(T)$, r\'eguli\`ere sur $k$, qui, par sp\'ecialisation (par exemple en $T=0$)
donne l'extension $K/k$?

Cette question, appel\'ee parfois probl\`eme de Beckmann-Black [$1$],
et que Black $[2]$ a conjectur\'ee avoir toujours une r\'eponse positive,
est r\'esolue, en particulier dans le cas o\`u $k=\Q$, pour les groupes
ab\'eliens, le groupe sym\'etrique, le groupe altern\'e $A_n$ et
une infinit\'e de groupes di\'edraux.

Nous montrons ici, comme corollaire du th\'eor\`eme $1$, qu'il en est de m\^eme pour le groupe $\PSL_2(\F_7)$.
Ce groupe est d'ordre $168$; il est isomorphe \`a $G=\PGL_3(\F_2)$,
et admet donc une repr\'esentation de degr\'e $7$. 

\ni
Plus pr\'ecis\'ement, nous montrons :

\medskip

{\sc Th\'eor\`eme 2.- } {\it Soit $(x_a)_{a\in A}$ sept ind\'etermin\'ees indic\'ees par
$A$, le plan projectif sur $\F_2$,
$K=k((x_a)_{a\in A})$  
et $P=\prod_a  (X-x_a)$.
Il existe un polyn\^ome $Q\in K[X]$, de degr\'e $6$, dont les coefficients
sont invariants par $G$, et 
tel que le groupe de Galois de $P-TQ$ sur le corps $K(T)$, o\`u $T$ est une nouvelle
ind\'etermin\'ee, est isomorphe \`a $G$. Le nombre de points de ramification du
rev\^etement $X\ms T=P/Q$ est \'egal \`a $6$, et le type de ramification 
est $(2^2,2^2,2^2,2^2,2^2,2^2)$.}

\medskip

d'o\`u l'on d\'eduit ais\'ement :

{\sc Th\'eor\`eme 2'.-} {\it Soit $k$ un corps de caract\'eristique $0$. 
Il existe un polyn\^ome 
non nul $H_1\in \Z[A_0,\ldots,A_6]$, o\`u les $A_i$ sont des ind\'etermin\'ees, 
tel que, si $P=X^7+a_6X^6+\ldots a_0$
est un \'el\'ement de $k[X]$ de groupe de 
Galois
contenu dans $G$ et  tel que $H_1(a_0,\ldots,a_6)\neq 0$, 
il existe $Q\in k[X]$ de degr\'e $\leq 6$, tel que
l'extension $L$ de $k(T)$ obtenue par adjonction des racines de $P-TQ$ soit
$k$-r\'eguli\`ere de groupe de Galois $G$.}

\medskip
{\bf Remarque.-} Il existe de m\^eme un polyn\^ome non nul 
$H_2\in \Z[A_0,\ldots,A_6]$
tel que, si les coefficients de $P$ n'annulent pas $H_2$, le type de ramification
de l'extension $L$ soit $(2^2,2^2,2^2,2^2,2^2,2^2)$.

\medskip
{\bf Exemples avec $k=\Q$.-}

1) On peut prendre pour $P$ un polyn\^ome scind\'e, par exemple
$P=X(X^2-1)(X^2-4)(X^2-9)$. On trouve alors
$Q=-259x^6-18578x^5+28812x^4+93494x^3-22709x^2-192228x-151380$.

2) Si l'on prend le polyn\^ome de Trinck $P=X^7-7X+3$, bien connu pour avoir 
$G$ comme groupe de Galois,
on trouve $$Q=2X^5-X^4-X^3-X^2-X+2=(X-1)^2(X+1)(2X^2+X+2).$$

Le discriminant de $P-TQ$ est $$81(T^2-5T+7)^2(800T^3+21T^2-441T+3087)^2.$$

3) On retrouve comme cas particuliers les familles de polyn\^omes 
de groupe de Galois 
$G$ trouv\'ees par
La Macchia $[3]$, ainsi que par Matzat et Malle $[4]$.

4) Dans $[5]$, Malle construit une famille de rev\^etements de $P^1$ dans $P^1$  de groupe de Galois $G$ ayant la m\^eme
dimension "g\'eom\'etrique" (\`a savoir $3$) que ceux construits ici.

\bigskip

\begin{center}
II. {\sc Probl\`eme de Noether}
\end{center}

\medskip

Comme autre corollaire du th\'eor\`eme $1$, nous montrons que le probl\`eme de
Noether admet une r\'eponse positive pour le groupe $G=\PGL_3(\F_2)$ agissant
sur sept points. Plus
pr\'ecis\'ement :

{\sc Th\'eor\`eme 3.-} {\it Soit $(x_a)_{a\in P^2(\F_2)}$ sept ind\'etermin\'ees indic\'ees par le plan
projectif sur $\F_2$; soit $K=\Q((x_a))$, muni de son action naturelle par $G$.
Le corps $K^G$ des points fixes de $K$ sous l'action de $G$ est une extension transcendante
pure de degr\'e $7$ de $\Q$.}

Comme il est bien connu, ceci est \'equivalent \`a l'existence d'un polyn\^ome 
$P\in \Q(T_1,\ldots,T_7)[X]$, o\`u les $T_i$ sont des ind\'etermin\'ees, tel que tout polyn\^ome
"suffisamment g\'en\'eral" de degr\'e $7$ d\'efini sur une extension
$k$ de $\Q$ et dont le groupe de Galois est contenu dans $G$ s'obtient en sp\'ecialisant les
ind\'etermin\'ees $T_i$ en des valeurs de $k$.

\medskip

\ni
{\sc Remarque.- } Dans un preprint $[7]$, B. Plants montre, en utilisant 
l'article $[6]$, o\`u est prouv\'e un analogue du th\'eor\`eme $2$ pour les
groupes altern\'es d'ordre impair, que le groupe
altern\'e $A_{2n+1}$ v\'erifie la propri\'et\'e de Noether  si et seulement si il en est de m\^eme de
$A_{2n}$. Il retrouve ainsi le fait que $A_5$ poss\`ede la propri\'et\'e de Noether.

\medskip

\ni

\section{Le plan projectif $A$}

Afin de pouvoir donner des formules explicites, num\'erotons les sept \'el\'ements 
$(1,\ldots,7)$ de $A$ et les sept \'el\'ements $(1',\ldots,7')$ 
de $B$ de la 
fa\c con suivante :
$$\begin{array}{ll}
1'=(234)\\
2'=(135)\\
3'=(126)\\
4'=(147)\\
5'=(257)\\
6'=(367)\\
7'=(456)
\end{array}$$

\ni
Ainsi, $a\in b'$ si et seulement si $b\in a'$.

\ni
Pour tout \'el\'ement $i$ de $A$ ou $B$
notons $G_i$ le sous-groupe de $G=\PGL_3(\F_2)$
stabilisant $i$;  
$G_i$ est isomorphe au groupe sym\'etrique $S_4$.
Par exemple, $G_1$ est engendr\'e par les deux permutations
$(234)(576)$ et $(2365)(47)$ et $G_{1'}$ par
$(234)(576)$ et $(1567)(34)$.

\ni
Soient $(x_i)_{1\leq i\leq 7}$, $X$ et $Y$  
des ind\'etermin\'ees,  
$P=\prod (X-x_i)$
et $K=k((x_i))$;  on fait agir $G$ 
sur $K(X,Y)$ trivialement sur $k(X,Y)$ et, pour tout $\sigma\in G$, par
$\sigma(x_i)=x_{\sigma(i)}, 1\leq i\leq 7$.

\ni
On choisit d'autre part un $7$-Sylow $Syl$ de $G$; pour tout $i\in A$ ou $B$,
l'application $Syl\ra G/G_i$ est bijective.

\ni
Nous aurons besoin
du lemme suivant, o\`u $\eps_i:\;G_i\ra \{\pm 1\}$ est la signature:

\ni
{\sc Lemme.-} {\it Pour tout $i$ dans $A$ ou $B$, notons 
$St_i$ le $k$-sous-espace 
vectoriel de
$k[x_1,\ldots,x_7]$ 
form\'e des polyn\^omes homog\`enes $p\in k[x_1,\ldots,x_7]$ de degr\'e 
au plus $1$ en chaque variable et de degr\'e total $3$, 
tels que, pour tout $\sigma\in G_i$, on ait
$\sigma(p)=\eps_i(\sigma) p$; $St_i$ est de dimension $1$. 
Plus particuli\`erement, le polyn\^ome $$\begin{array}{r}
u_1=(x_2-x_3)(x_4-x_5)(x_6-x_7)+(x_2-x_4)(x_3-x_7)(x_5-x_6)
\end{array}$$ 
est un g\'en\'erateur de $St_1$,
et
$$v_1=x_2(x_5-x_7)(x_1-x_6)+x_3(x_1-x_5)(x_7-x_6)+x_4(x_1-x_7)(x_6-x_5)$$
est un g\'en\'erateur de $St_{1'}$.
De plus, pour tout $a \in k$, 
on a $v_1(x_1+a,\ldots,x_7+a)=v_1(x_1,\ldots,x_7)$.
}

{\sc Remarques.-} 1) $u_1$ est nul (pour des sp\'ecialisations $\alpha_i$ distinctes deux-\`a-deux des ind\'etermin\'ees $x_i$) si et seulement s'il existe une homographie involutive
permutant $\alpha_2$ et $\alpha_6$, $\alpha_3$ et $\alpha_5$ , $\alpha_4$ et $\alpha_7$.

2) Par suite, il n'est pas possible, pour des sp\'ecialisations $\alpha_i$ des $x_i$ distinctes deux-\`a-deux,
que  $u_1$ et ses analogues $u_i$ soient tous nuls : sinon, on aurait sept homographies involutives commutant  
deux-\`a-deux.

\section{D\'emonstration du th\'eor\`eme $1$}

Soit $F$ une correspondance $3-3$ compatible avec les relations droites-plans de $P^2(\F2)$,
c'est-\`a-dire un polyn\^ome \`a deux variables $X,Y$ de bi-degr\'e $3$, d\'efini \`a une constante
multiplicative pr\`es.
On peut, par composition de la correspondance $F$ avec elle-m\^eme,  lui 
associer une correspondance ``compl\'ementaire'' $H$, de bidegr\'e $4-4$,
telle que, si $(x_a)_{a\in A}$ et $(y_b)_{b\in B}$ sont des quantit\'es telles que $F(x_a,y_b)=0$ si et seulement
si $a\in b$, on ait $H(x_a,y_b)=0$ si et seulement si $a\not\in b$. 

Plus pr\'ecis\'ement,  le r\'esultant de $F(X,Y)$ et $F(Z,Y)$ relativement \`a $Y$ est de la forme
$(X-Z)^3R_1(X,Z)$, o\`u
$R_1$ est un polyn\^ome sym\'etrique de bidegr\'e $6-6$, qui associe \`a l'un des $x_a$ les six autres.
De m\^eme, le r\'esultant de $R_1(X,Z)$ et de $F(Z,Y)$ relativement \`a $Z$ est de la forme
$F(X,Y)^2H(X,Y)^3$, o\`u $H$ est la correspondance cherch\'ee.

Pour montrer le th\'eor\`eme $1$, nous prouvons d'abord l'existence de $H$, puis
celle de $F$ :

\medskip

{\sc Th\'eor\`eme 1''.-} {\it 
Conservons les notations du th\'eor\`eme $1$. Il existe une famille
$(y_j)_{j\in B}$ d'\'el\'ements de $K$,
une $3-3$ correspondance $\Gamma$, 
et une $4-4$ correspondance $\Gamma'$  sur $P^1\times P^1$,
telle que $(x_i,y_j)\in \Gamma$  (resp. $(x_i,y_j)\in \Gamma'$)
si et seulement si $i\in j$ (resp. $i\not\in j$).}

Le th\'eor\`eme $1''$ revient \`a montrer 
l'existence de deux polyn\^omes $F$ et $H$ dans $K[X,Y]$, 
de
bidegr\'es respectifs $3-3$ et $4-4$, et de sept \'el\'ements $(y_1,\ldots,y_7)$ de $K$, tels que
$F(x_i,y_j)=0$ si et seulement si $i\in j'$ et $H(x_i,y_j)=0$ 
si et seulement si
$i\not\in j'$.

\ni
Pour tout $j$, $1\leq j\leq 7$, notons 
$r_j=\prod_{i\in j'} (X-x_i)$ et $s_j=\prod_{i\not\in j'}(X-x_i).$

\ni
Les assertions du th\'eor\`eme reviennent \`a dire qu'il existe
trois suites $(y_j)$, $(a_j)$ et $(b_j)$, $1\leq j\leq 7$, telles que,
si l'on pose, pour $1\leq j\leq 7$, $L_j(Y)=\prod_{k\neq j} (Y-y_k)$, 
le polyn\^ome 
$$F=\sum a_j L_j(Y) r_j(X),\;\;\; {\rm resp.}\;\;\; 
H=\sum b_j L_j(Y)s_j(X)$$ 
est de degr\'e $3$ (resp. $4$) en $Y$.

\ni
Le fait que $H$ est de degr\'e au plus $4$ en $Y$ s'\'ecrit
$$\left\{\begin{array}{l}
\sum b_j s_j=0\\
\sum b_j y_j s_j=0\end{array}\right.$$

\ni
Comme les sept polyn\^omes $s_j$ sont de degr\'e $4$, l'espace des
combinaisons lin\'eaires les annulant est de dimension  $2$ .
En effet, si l'on choisit cinq $s_j$, il est clair que quatre d'entre eux exactement  ont une racine commune,
et le cinqui\`eme ne peut donc pas en \^etre une combinaison lin\'eaire. 
Supposons pour fixer les id\'ees que $x_1$ est cette racine commune; les quatre polyn\^omes sont alors
$s_1,s_5, s_6$ et $s_7$. Le d\'eterminant de ces quatre polyn\^omes divis\'es par $X-x_1$
dans la base canonique $\{1,X,X^2,X^3\}$ est  \'egal, au signe pr\`es, \`a
$$u_1(x_2-x_6)(x_3-x_5)(x_4-x_7),$$
o\`u $u_1$ a \'et\'e d\'efini dans la section pr\'ec\'edente, et est donc non nul. Par suite, toute partie
\`a cinq \'el\'ements de $\{s_1,\ldots,s_7\}$ est libre.

 Soit donc deux combinaisons lin\'eaires ind\'ependantes
$\sum_j b_j s_j=\sum_j b'_j s_j=0,$
et posons $y_j=b'_j/b_j$ pour $1\leq j\leq 7$.

Les $y_i$ sont distincts deux-\`a-deux (car si $y_i$ \'etait \'egal \`a $y_j$, la famille des
$(s_k)$, $k\neq i,j$ serait li\'ee). Si $\{(c_j),(c'_j)\}$ est une autre base des combinaisons
lin\'eaires reliant les $(s_i)$, la matrice de passage de la base $\{b_j),(b'_j)\}$ vers la pr\'ec\'edente
donne l'homographie reliant les $(c'_j/c_j)$ aux $y_j$.

\ni
Une fois les $y_j$  trouv\'es,  la d\'etermination
des $(a_j)$ revient \`a r\'esoudre le syst\`eme lin\'eaire
\`a $7$ inconnues et $12$   \'equations
$$\sum a_j r_j=\sum a_j y_j r_j=\sum a_j y_j^2 r_j=0,$$
qui s'av\`ere 
avoir une unique solution (\`a un scalaire multiplicatif pr\`es).

\bigskip

\ni
Remarquons que 
$F,H$ et les $y_i$ ne sont pas uniques, puisqu'\`a partir
de polyn\^omes $F,H$ et d'\'el\'ements $y_i$ satisfaisant  aux relations
d'incidence,
on peut en obtenir   une infinit\'e d'autres
en les transformant par une homographie $Y\ms (aY+b)/(cY+d)$, $a,b,c,d\in K$. 

Cependant, une fois fix\'es trois $y_i$, $F$ est unique.
On peut par exemple prendre 
$y_1=\infty,y_2=0,y3=1$. N\'eanmoins, dans ce cas, les coefficients de $F$ et $H$
ne sont pas stables par l'action de $G$. 

Si, par contre, on peut trouver des $y_i\in K$ tels que, pour tout $\sigma\in G$, $\sigma(y_k)=
y_{\sigma(k)}$, les coefficients de $F$ (d\'efinis \`a une constante pr\`es) peuvent
\^etre choisis dans $K^G$, le sous-corps de $K$ fix\'e par $G$.

Par ailleurs, si deux solutions $(y_i)$ et $(z_i)$ sont telles que la relation ci-dessus soit v\'erifi\'ee
pour chacune d'elles, l'homographie permettant de passer de l'une \`a l'autre
peut aussi \^etre choisie \`a coefficients dans $K^G$.

\medskip

Nous montrons ci-apr\`es qu'il existe une famille
 famille $(y_i)$ plus satisfaisante
que les autres, du point de vue galoisien et de la compatibilit\'e avec les homographies sur les
$x_i$ : 

\ni
Posons 
$$\left\{\begin{array}{l}
y_1=-\fracb{v_1(1/x_1,\ldots,1/x_7) x_1\ldots x_7}{v_1}, \\
y_{\sigma(1)}=\sigma(y_1) \;\;{\rm et}\;\;u_{\sigma(1)}=\sigma(u_1)\;\;
{\rm pour}\;\;\sigma \in Syl\\
U=\prod_{i=1}^7 (Y-y_i)\\  
l_1=u_2u_3u_4v_1^2\\  
m_1=u_2u_3u_4v_1(x_2-x_3)(x_3-x_4)(x_4-x_2)
\end{array}\right.
.$$ 

{\sc Proposition 1.-} {\it 
Les $y_i$ construits ci-dessus v\'erifient les propri\'et\'es 
suivantes:

a) pour tout $\sigma\in G$, $y_{\sigma(i)}=\sigma(y_i)$

b) Pour toute homographie $h:\;z\ms (az+b)/(cz+d)$, $a,b,c,d\in k$,
on a
$$y_i(h(x_1),\ldots,h(x_7))=h(y_i(x_1,\ldots,x_7)).$$

c) l'application $(x_1,\ldots,x_7)\ms (y_1,\ldots,y_7)$ est involutive.

Le th\'eor\`eme $1'$ est alors v\'erifi\'e 
en prenant pour $\Gamma$ (resp. $\Gamma'$)  la correspondance
d'\'equation $F=0$ (resp. $H=0$), avec
$$F=U\sum _{\sigma\in Syl} \sigma(\fracb{l_1 r_1}{Y-y_1})\;\;{\rm et}\;\;
H=U\sum_{\sigma\in Syl} \sigma(\fracb{m_1 s_1}{Y-y_1}).$$

De plus, les coefficients de $F$ et $H$  
sont invariants par $G$.
}

\bigskip
En effectuant la division
euclidienne de $FH$ (vu comme polyn\^ome en $X$) par $P$, on
obtient $FH=P(X)V(Y)-\phi(X,Y)$, o\`u $\phi$ est de degr\'e $\leq 6$ en
$X$ et de degr\'e $\leq 7$ en $Y$. Comme, pour toute racine $y$ de $U$,
$\phi(X,y)$ s'annule en les $7$ racines de $P$, $\phi$ est de la forme
$Q(X)U(Y)$, $Q$ de degr\'e $\leq 6$, et on a
$$FH=P(X)V(Y)-Q(X)U(Y).$$

Par ailleurs, on a la formule explicite suivante pour $Q$ :

{\sc Proposition 2.-} {\it 
On  a 
$$Q=P\sum_{\sigma\in Syl} \sigma(\frac{q_1}{X-x_1}),$$
o\`u $q_1=u_1^2 P'(x_1)(x_2-x_6)(x_3-x_5)(x_4-x_7)$.}

Il est alors clair que les coefficients de $Q$ sont invariants par $G$, et que $P$ et $Q$ sont premiers
entre eux.

\medskip

Si $T$ est une ind\'etermin\'ee, et si $P_T=P-TQ$ et $U_T=U-TV$, on
a imm\'ediatement $FH=P_T(X)V(Y)-Q(X)U_T(Y).$

On en d\'eduit que l'ensemble $A$ des racines de $P_T$ et l'ensemble $A'$ des
racines de $U_T$ sont reli\'ees par une relation d'incidence de
type ``points-droites'' sur le plan projectif $P^2(\F_2)$ : si $y\in A'$, les trois ``points''
de $y$ sont les racines de $F(X,y)$. 

Ceci prouve le th\'eor\`eme $1$.

Le polyn\^ome $S$ du th\'eor\`eme $1'$ s'obtient par exemple en multipliant le discriminant de $P$ par
le produit des $\sigma(u_1)$, o\`u $\sigma$ parcourt le groupe sym\'etrique $S_7$, et en exprimant
le r\'esultat en fonction des coefficients de $P$.

\medskip

{\sc Remarques.-} 
1) Soit $disc(P,X)$ le discriminant d'un polyn\^ome  relativement \`a $X$.
Apr\`es calculs, il s'av\`ere que $$
(\prod_b v_b)^{12} disc(P-TQ,X)=(\prod_a u_a)^{12} disc (U-TV,Y).$$

Cela prouve que,
si $P\in k[X]$ n'annule pas $S$, il en est de m\^eme de tout polyn\^ome  s\'eparable du faisceau
$P-TQ$ associ\'e .

 2) Si $F$ est la correspondance de la proposition $1$, 
on peut montrer que, pour toute sp\'ecialisation $t$ de $T$,
les racines $(\alpha_a)_{a\in A}$ de $P-tQ$ et $(\beta_b)_{b\in b}$ de $U-tV$
sont
telles que, pour tout $j\in B$, $\beta_j$ est la sp\'ecialisation en les 
$\alpha_i$ de la fraction rationnelle 
$y_j(x_1,\ldots,x_7)$.  

\section{D\'emonstration des th\'eor\`emes $2$ et $2'$}

Le groupe de Galois de  l'extension $M=K(T)(A)$ de $K(T)$ est contenu
dans $G$ : si $\sigma$ est un \'el\'ement de ce groupe,
et si trois \'el\'ements de $A$ sont ``align\'es'', i.e. s'il
existe $y\in A'$ tel que ces trois \'el\'ements sont les racines
de $F(X,y)$, leurs images par $\sigma$ sont les racines
de $F(X,\sigma(y))$, et sont donc align\'es; $\sigma$ pr\'eservant
l'alignement, il est dans $G$.

Par ailleurs, on peut montrer , par exemple en prenant l'un des exemples donn\'es en introduction, que le type de ramification de
$X\ms P(X)/Q(X)$ est $(2^2,2^2,2^2,2^2,2^2,2^2)$. 
Comme le seul sous-groupe transitif de $G\subset S_7$ engendr\'e par
des produits de deux transpositions est $G$ lui-m\^eme, on
en d\'eduit que
$G$ est le groupe de
Galois de $M/K(T)$, d'o\`u le th\'eor\`eme $2$.

Le th\'eor\`eme $2'$   s'en d\'eduit facilement : si un polyn\^ome
s\'eparable
$P$ a un groupe de Galois contenu dans $G$, on peut mettre en
bijection l'ensemble de ses racines et le plan projectif $A$,
de fa\c con que les coefficients de $F$ , $G$ et $Q$ appartiennent
au corps de base $k$. Par un raisonnement analogue au pr\'ec\'edent,
on montre alors que $G$ est le groupe de Galois de $P-TQ$ sur $k(T)$,
pourvu que le r\'esultant $R$ de $P$ et $Q$ soit non nul; on peut alors
prendre $H_1=\prod_{\sigma\in S_7/G} R$. 
Le discriminant de $P-TQ$ relativement \`a $X$ est le carr\'e
d'un polyn\^ome $S(T)$ de degr\'e $6$, dont le discriminant 
permet d'obtenir le polyn\^ome $H_2$ de la remarque
suivant le th\'eor\`eme $3$.

Remarquons d'autre part que, g\'en\'eriquement, le polyn\^ome $Q$ a un discriminant non nul
(comme le montre sa sp\'ecialisation en le premier exemple de l'introduction). En particulier,
d'apr\`es la remarque $1$ de la section pr\'ec\'edente, le polyn\^ome $Q$ unitaire construit \`a partir du polyn\^ome $\prod_{a\in A} (X-x_a)$, o\`u les
$x_a$ sont des ind\'etermin\'ees, v\'erifie les hypoth\`eses du
th\'eor\`eme $1'$. 

\section{Le probl\`eme de Noether pour $\L_3(2)$}

Le stabilisateur d'un point  $a\in A$ dans $G=\L_3(2)$ est
isomorphe au groupe sym\'etrique $S_4$, agissant transitivement sur les six parties \`a deux
\'el\'ements de $\{1,2,3,4\}$. Il est bien connu que le probl\`eme de
Noether a une solution positive dans ce cas, c'est-\`-dire qu'il existe
un polyn\^ome g\'en\'erique unitaire ${\cal Q}$, \`a six ind\'etermin\'ees $T_1,\ldots,T_6$, tel que tout
polyn\^ome unitaire de degr\'e $6$ suffisamment g\'en\'eral sur un corps $k$ de caract\'eristique nulle, de groupe de Galois contenu dans $S_4$
via cette action, 
est obtenu par sp\'ecialisation de ces
ind\'etermin\'ees en des \'el\'ements de $k$. Une fa\c con d'obtenir ${\cal Q}$ est par exemple la suivante : si $P\in k[X]$ de degr\'e $6$ a comme groupe de Galois $S_4\subset A_6$, ses racines
se regroupent en trois couples $(x_1,x'_1), (x_2,x'_2)$ et $(x_3,x'_3)$ telles que
$x'_i\in k(x_i)$. Le polyn\^ome dont les trois racines sont les
$z_i=(x_i-x'_i)^2$ est  de la forme $X^3-a_1X^2+a_2X-u^2$, $a_1,b_1$ et $u\in k$,
et, comme $x_i+x'_i\in k(z_i^2)$,  il existe $a_3,a_4$ et $a_5$ dans $k$ tels que $x_i=z_i+a_3+a_4 z_i^2+a_5 z_i^4;$ Le polyn\^ome g\'en\'erique est le polyn\^ome minimal de $z$.
 Les quantit\'es $a_1,a_2,a_3,a_4,a_5$ et $u$, exprim\'ees en fonction
des racines d'un polyn\^ome ind\'etermin\'e $P=(X-x_1)...(X-x_6)$, engendrent le sous-corps
de $\Q(x_1,\ldots,x_6)$ fix\'e par $S_4$.

Soit donc $(x_a)_{a\in A}$ sept ind\'etermin\'ees indic\'ees par $A=P^2(\F_2)$,
$K=\Q((X_a))$,  $M=K^G$, et soit $F$ la correspondance d\'efinie sur $M$
de la proposition $1$. 
Soit $Q_1\in M[X]$ le polyn\^ome unitaire proportionnel \`a $Q$;
par construction, il existe une bijection
$a\ra \alpha_a$ de $A$ sur l'ensemble form\'e des 
six racines de $Q_1$ et de l'infini et une bijection $b\ra (\beta_b)$ 
de $B$ sur les racines de $V$
telles que $F(\alpha_a,\beta_b)=0$ si et seulement si $a\in b$. Le
groupe de Galois de $Q_1$ sur $M$ fixe l'infini et permute les autres 
$(\alpha_a)$, il est donc contenu dans le stabilisateur $St_{\infty}$ de l'infini dans $\L_3(2)$, et
lui est en fait \'egal, puisque c'est le cas pour l'exemple $1$ de 
l'introduction, qui en est une sp\'ecialisation.

Par ailleurs, comme on l'a vu pr\'ec\'edemment, le 
polyn\^ome $Q_1$ est suffisamment g\'en\'eral, dans la terminologie du
th\'eor\`eme $1'$, et la correspondance $F$ pr\'eservant les relations
d'incidence entre les $(\alpha_a)$ et les $(\beta_b)$ est donc unique.

Par ailleurs, $Q_1$ s'obtient par sp\'ecialisation de $T_1,\ldots,T_6$ en des \'el\'ements
$r_1,\ldots,r_6\in M$ du polyn\^ome g\'en\'erique  ${\cal Q}$;
le polyn\^ome $Q_1$ est donc \`a coefficients dans
$L=\Q(r_1,\ldots,r_6)$, a $St_{\infty}$ comme groupe de Galois sur $L$.
D'apr\`es la remarque $2$ de la fin de la section $3$,
$F$ est donc donn\'ee par la formule de la proposition $1$, o\`u les $x_i$
sont sp\'ecialis\'es en les $\alpha_i$, 
et $F$ est \`a coefficients dans $L$.

Le polyn\^ome initial
 $P=\prod_{a\in A}(X-X_a)$ s'obtient \`a partir de $F$ en \'eliminant $Y$ entre les
 relations $F(X_a,Y)$ et $F(X,Y)$, ($a\in A$), et est donc \`a coefficients dans $L(X_a)$, pour tout $a\in A$; le groupe de Galois
 de $K/L(X_a)$ est \'egal au stabilisateur de $a$ dans $G$, et donc le groupe de Galois
 de l'intersection $N=\cap_{a\in A}L(X_a)$ est \'egal \`a $G$, d'o\`u $N=K^G$.
 Comme, pour $a\in A$, $N$ est une sous-extension de l'extension $L(X_a)/L$, qui est transcendante pure de 
degr\'e $1$,
il existe $r_7\in N$ tel que $r_1,\ldots,r_7$ sont alg\'ebriquement ind\'ependants, 
et telle que $K^G=N=L(r_7)=\Q(r_1,\ldots,r_7)$, d'o\`u le r\'esultat.
 
De fa\c con explicite, on peut prendre $r_7=\prod_{a\in A} x_a$ : soit en effet $R\in L[X]$ le polyn\^ome unitaire de degr\'e $7$,
 s'annulant en $0$,  associ\'e \`a $F$. 
Il existe $t\in M$ tel que le polyn\^ome $P=\prod (X-x_a)$ s'\'ecrit
 $R-tQ_1$. Comme $R$ et $Q_1$ sont unitaires, on a $-\prod x_a=R(0)-tQ_1(0)=-tq_6$,
 o\`u $q_6$ est 
 le terme constant de $Q_1$. Par suite, $t= \frac{\prod x_a }{q_6}$, d'o\`u le r\'esultat. 
 Par construction, le polyn\^ome $R-T_7Q_1 \in \Q(T_1,\ldots,T_7)[X]$ est un polyn\^ome g\'en\'erique
 pour les polyn\^omes de degr\'e $7$ dont le groupe de Galois est un sous-groupe de $\L_3(2)$.
 
\section{D'autres cas}
Pour terminer, donnons 
quelques r\'esultats analogues aux th\'eor\`emes $1$ et $2$ 
pour les polyn\^omes de degr\'e
$4$, $5$ et $6$:

\medskip

1) {\it Soit $P\in k[X]$ un polyn\^ome de degr\'e
$4$, suffisamment g\'en\'eral (en ce sens que ses coefficients n'annulent
pas un certain polyn\^ome non nul),  
et dont le groupe de Galois est contenu dans $(\Z/2\Z)^2$. 
Il existe un polyn\^ome $Q\in k[X]$ de degr\'e $4$
tel que le groupe de Galois de $P-TQ$ sur $k(T)$ est $(\Z/2\Z)^2$.}

\medskip

Ici, le polyn\^ome $Q$ est particuli\`erement simple \`a calculer:
si on homog\'en\'eise $P$ par $R(X,Z)=Z^4P(X/Z)$, 

et si $$Hes=R''_{X^2}R''_{Z^2}-(R''_{XZ})^2,$$ on a $Q=Hes(X,1).$
\`A une constante multiplicative pr\`es, $Q$ est le reste de $P'^2\mod P.$

Par ailleurs, $P(X)Q(Y)-P(Y)Q(X)=(X-Y)F_1(X,Y)F_2(X,Y)F_3(X,Y),$
o\`u $F_1$, $F_2$ et $F_3$ sont des correspondances $1-1$ sym\'etriques,
en fait les graphes des trois homographies involutives envoyant une racine
de $P-TQ$ sur ses trois autres conjugu\'ees. Ces homographies sont
ind\'ependantes de $T$.  

\bigskip

 2) {\it Soit $P\in k[X]$ un polyn\^ome de degr\'e
$4$, suffisamment g\'en\'eral 
et dont le groupe de Galois est contenu dans le groupe di\'edral $D_4$.
Il existe deux polyn\^omes $Q1$ et $Q_2 \in k[X]$ de degr\'e $\leq 4$

tel que le groupe de Galois de $P-T_1Q_1-T_2Q_2$ sur $k(T)$ soit $D_4$.}

Explicitement, si $P=\prod_{i \in \Z/4\Z} (X-x_i)$,  si $u=x_0-x_1+x_2-x_3$ 
et si
$q=ux+x_1 x_3-x_0x_2$,
on peut prendre
$Q_1=q^2$ et $Q_2=uq (2x^2-x(x_0+x_1+x_2+x_3)+x_0x_2+x_1x_3),$

$F$ \'etant une correspondance (\`a un param\`etre) $2-2$ telle que $F(x_i,y_j)=0$
si et seulement si $j=i$ ou $j=i+1$, et $H$ la correspondance
 $1-1$ graphe de l'homographie involutive 
envoyant, pour tout $i$, $x_i$ sur $x_{i+2}$.  

\medskip

3) 
{\it Soit $P\in k[X]$ de degr\'e $5$, suffisamment g\'en\'eral, tel que
le groupe de Galois de $P$ soit contenu dans le groupe di\'edral $D_5$.
Il existe $Q\in k[X]$ de degr\'e $\leq 4$ tel que le groupe de Galois de $P-TQ$ sur $k(T)$
soit \'egal \`a $D_5$.}

\medskip

Soit $P=\prod_{i \in \Z/5\Z}(X-x_i).$ 
On a alors $$Q=P\sum_{i=0}^4 \frac{P'(x_i)(x_{i-1}-x_{i+1})(x_{i-2}-x_{i+2})}{X-x_i},$$
et $P(X)Q(Y)-P(Y)Q(X)=(X-Y)FG,$
o\`u $F$ et $G$ sont les deux correspondances sym\'etriques
de bidegr\'e $2-2$ telles que $F(x_i,x_j)=0$ (resp. $G(x_i,y_j)=0$) 
si et seulement si 
$i-j=\pm 1$ (resp. $i-j=\pm 2.$) 

\medskip

4) 
{\it Soit $P\in k[X]$ de degr\'e $6$, suffisamment g\'en\'eral, tel que
le groupe de Galois de $P$ soit contenu dans le  produit en couronnes 
$G$ de $S_3^2$ par $\Z/2\Z$ (d'ordre $72$).
Il existe $Q_1$ et $Q_2\in k[X]$ de degr\'e $5$ tel que le groupe de Galois 
de $P-T_1Q_1-T_2Q_2$ sur $k(T_1,T_2)$
soit \'egal \`a $G$.}

\medskip

{\sc Remarques.-} 1) Pour $n\geq 6$, il n'est pas vrai que, 
si $(x_i)_{i\in \Z/n\Z}$ sont des ind\'etermin\'ees, il existe une
$2-2$ correspondance $F$   compatible avec les relations d'incidence 
associ\'ees \`a un polygone r\'egulier \`a $n$ c\^ot\'es (dont le groupe associ\'e est le groupe di\'edral $D_n$) d\'efinies dans l'introduction: on sait d'apr\`es l'interpr\'etation
du th\'eor\`eme de Poncelet par Jacobi que ceci revient \`a se donner une
courbe elliptique munie d'un point d'ordre $n$ et d'un autre point,
situation param\'etr\'ee par 
une surface, alors qu'ici, une fois $x_1,x_2$ et $x_3$
ramen\'es \`a $\infty,0$ et $1$, on a une vari\'et\'e de dimension
$n-3$. Pour $n=5$, cela explique qu'il n'y a pas de condition sur les
$x_i$ (point $3)$ ci-dessus), et , pour $n=4$, qu'on trouve une famille \`a deux param\`etres au
point $2)$ ci-dessus.

Par contre, soit $(x_1,\ldots,x_6)$ six quantit\'es. Il n'existe en g\'en\'eral pas de correspondance
$F$ de bidegr\'e $2-2$ et six quantit\'es $(y_1,\ldots,y_6)$ telles que 
$F(x_i,y_j)=0$ si et seulement si $j=i$ ou $j=i+1$. Pour qu'il en existe, il faut et il suffit que
 
$$(x_0-x_1)(x_2-x_3)(x_4-x_5)+(x_1-x_2)(x_3-x_4)(x_5-x_0)=0,$$
i.e. qu'il existe une homographie involutive transformant
$x_i$ en $x_{i+3}$ pour $i\in \Z/6\Z$.

Si cette condition est r\'ealis\'ee, il existe $Q$ de degr\'e $\leq 5$
tel que $P-TQ$ ait $D_6$ comme groupe de Galois sur $\Q(x_0,\ldots,x_5)(T)$.

\medskip

2) Soit $F\in k[X,Y]$ une correspondance de bi-degr\'e $2-2$, associ\'ee \`a la configuration
de Poncelet associ\'ee au polygone \`a $n$ c\^ot\'es; $F$ est compatible si et seulement si
le diviseur $D_X-D_Y$ est d'ordre $n$, o\`u $D_X$ (resp. $D_Y$) est la partie polaire du diviseur
de la fonction $X=0$ (resp. $Y=0$). 

Soit $F$ une correspondance $3-3$ compatible avec $P^2(\F_2)$.  Si $F(u,v)=0$, le diviseur de la fonction $(Y-v)/(X-u)$ est
\'egal \`a $7(D_X-D_Y)$, donc  $D_X-D_Y$ est d'ordre $7$ dans le groupe de Picard de la courbe (g\'en\'eriquement de genre $4$)
d'\'equation $F=0$. Mais cette
condition n'est pas suffisante . Par exemple, la courbe 

$$\begin{array}{ll}F=&
{X}^{3}{Y}^{3}+2\,{X}^{3}{Y}^{2}+4\,{X}^{3}Y+4\,{X}^{3}+3\,{X}^{2}{Y}^
{3}+4\,{X}^{2}Y+\\&3\,X{Y}^{3}+2\,X{Y}^{2}+4\,XY+4\,X+2\,{Y}^{2}+4\,Y+2\end{array}$$
d\'efinie sur $\F_5$ est une correspondance $3-3$, le diviseur $D_X-D_Y$ est d'ordre $7$, mais elle
n'est pas $P^2(\F_2)$ compatible.

\bigskip

{\sc R\'ef\'erences.-} 

[1] S. Beckmann, Is every extension of $\Q$ the specialization of a branched
covering?, J. Algebra 164 (1994), 430-451.  

[2] E. Black, On semi-direct products and the arithmetic lifting property,
J. London Math. Soc., (2) , 60 , (1999), 677-688.

[3] S.E. LaMacchia, Polynomials with Galois group $\PSL(2, 7)$,
Communications in Algebra,8,(1980), 983-992.

[4] G. Malle and B.H. Matzat, Inverse Galois Theory, Springer-Verlag,
Berlin, Heidelberg, New York, 1999.

[5] G. Malle, Multi-parameter Polynomials with Given Galois Group,
J. Symbolic Comput. 30 (2000), 717-731.

[6] J.-F. Mestre,  Extensions r\'eguli\`eres de $\Q(T)$ de groupe de Galois $\tilde{A}_n$,              Journal of Alg. 131 (1990), 483-495.

[7] Bernat Plans , On the $\Q$-rationality of $\Q(X_1,\ldots,X_5)^{A_5}$, preprint.

\end{document}